\input amstex\documentstyle{amsppt}  
\pagewidth{12.5cm}\pageheight{19cm}\magnification\magstep1
\topmatter
\title On the representations of disconnected reductive groups over $F_q$
\endtitle
\rightheadtext{Representations of disconnected reductive groups over $F_q$}
\author G. Lusztig\endauthor
\address{Department of Mathematics, M.I.T., Cambridge, MA 02139}\endaddress
\thanks{Supported in part by the National Science Foundation}\endthanks
\endtopmatter   
\document
\define\bcb{\bar{\cb}}
\define\Irr{\text{\rm Irr}}

\define\rank{\text{\rm rank}}

\define\dl{\dot l}
\define\dsg{\dot{\sg}}
\define\dW{\dot W}
\define\dH{\dot H}
\define\dD{\dot D}
\define\dQ{\dot Q}
\define\dJ{\dot J}
\define\dT{\dot T}
\define\dR{\dot R}
\define\dC{\dot C}
\define\dP{\dot P}
\define\dI{\dot I}

\define\dw{\dot w}
\define\dr{\dot{\r}}
\define\dt{\dot t}

\define\dc{\dot c}

\define\dce{\dot{\ce}}
\define\dcc{\dot{\cc}}
\define\dcar{\dot{\car}}
\define\dfA{\dot{\fA}}

\define\uz{\un z}

\define\uhx{\un{\hx}}

\define\hx{\hat x}

\define\da{\dagger}
\define\dsv{\dashv}

\define\si{\sim}
\define\wt{\widetilde}
\define\sqc{\sqcup}

\define\qua{\quad}

\define\hG{\hat G}

\define\bG{\bar G}

\define\bX{\bar X}

\define\lb{\linebreak}

\define\op{\oplus}
   
\redefine\sp{\spadesuit}
\define\part{\partial}

\define\n{\notin}

\define\m{\mapsto}
\define\do{\dots}

\define\sub{\subset}    

\define\T{\times}
\define\ti{\tilde}
\define\nl{\newline}
\redefine\i{^{-1}}

\define\un{\underline}
\define\ov{\overline}
\define\ot{\otimes}
\define\bbq{\bar{\QQ}_l}

\define\Ad{\text{\rm Ad}}
\define\Hom{\text{\rm Hom}}

\define\sg{\text{\rm sgn}}
\define\tr{\text{\rm tr}}

\define\di{\diamond}

\define\a{\alpha}
\redefine\b{\beta}
\redefine\c{\chi}

\redefine\d{\delta}

\redefine\o{\omega}
\define\p{\pi}
\define\ph{\phi}

\define\r{\rho}
\define\s{\sigma}

\define\th{\theta}

\redefine\l{\lambda}

\define\x{\xi}

\redefine\D{\Delta}

\define\Ph{\Phi}

\define\kk{\bold k}

\define\CC{\bold C}

\define\HH{\bold H}

\define\NN{\bold N}

\define\QQ{\bold Q}

\define\ZZ{\bold Z}

\define\ca{\Cal A}
\define\cb{\Cal B}
\define\cc{\Cal C}

\define\ce{\Cal E}
\define\cf{\Cal F}

\define\ch{\Cal H}
\define\ci{\Cal I}

\define\cl{\Cal L}

\define\co{\Cal O}

\define\car{\Cal R}

\define\fc{\frak c}

\define\fA{\frak A}

\define\tK{\ti K}
\define\tL{\ti L}

\define\tP{\ti P}

\define\DL{DL1}
\define\DLL{DL2}
\define\DLLL{DL3}
\define\DM{DM}
\define\KL{KL1}
\define\KLL{KL2}
\define\LC{L1}
\define\ORA{L2}
\define\LU{L3}
\define\MA{Ma}
\define\SPA{Sp}
\head Introduction\endhead
Let $G$ be a connected reductive algebraic group defined over a finite field $F_q$.
One of the main tools in the study of representations of the finite group $G(F_q)$ over a field of 
characteristic zero is the use of certain varieties $X_w$ (see \cite{\DL}) on which $G(F_q)$ acts (here $w$ is a
Weyl group element). Now let $\s:G@>>>G$ be a quasisemisimple automorphism of $G$ and let $m\ge1$ be an integer
such that $\s^m=1$. Consider the semidirect product $\hG$ of $G$ with the cyclic group of order $m$ with
generator $\s$; this is naturally an algebraic group defined over $F_q$. Now the finite group $\hG(F_q)$ acts 
naturally on the disjoint union $\sqc_w X_w$ and from this one can again derive information about the
representations of $\hG(F_q)$. For example, this observation has been used by the author in his proof of the 
finiteness of the number of unipotent $G$-conjugacy classes in the connected component $G\s$ of $\hG$ (see 
\cite{\SPA, I,4.1}); the connection between the varieties $X_w$ and the representations of $\hG(F_q)$ has been 
systematically investigated by Digne and Michel \cite{\DM} and by Malle \cite{\MA}. In this paper we try to 
extend some results on unipotent representations established for $G(F_q)$ in \cite{\ORA} to $\hG(F_q)$. 
One of the key steps in the description \cite{\ORA} of the set of unipotent representations of $G(F_q)$ is the
definition of a partition of that set into subsets indexed by the two-sided cells \cite{\KL}
of the Weyl group such that certain explicit $\QQ$-linear combinations of virtual 
representations of $G(F_q)$ given by the alternating sum of the cohomologies of the various $X_w$ are 
linear combinations of unipotent representations corresponding to a fixed two-sided cell. A conjectural 
extension of this statement to the case of $\hG(F_q)$ was formulated by Malle in \cite{\MA} and proved in this 
paper (see 2.4(ii)). Our proof is a generalization of that in \cite{\ORA}; the main new ingredient is the use of 
the generalization given in \cite{\LC} to certain Hecke algebras with unequal parameter of the polynomials 
$P_{y,w}$ defined in \cite{\KL} and of their geometric interpretation stated in \cite{\LC} (and generalizing 
that in \cite{\KLL}) which is proved in this paper. It turns out that these generalized polynomials appear 
naturally in the study of the varieties $X_w$ in connection with a nontrivial $G(F_q)$-coset of $\hG(F_q)$ and 
they provide the necessary tools to prove the above conjecture.

{\it Notation.}   
Let $\kk$ be an algebraic closure of the finite field $F_p$ with $p$ elements. If $q$ is a power of $p$ we denote
by $F_q$ the subfield of $\kk$ with $q$ elements. Let $\bbq$ be an algebraic closure of the field of $l$-adic 
numbers ($l$ is a fixed prime number $\ne p$). All algebraic varieties in this paper are over $\kk$. For an 
algebraic variety $Y$ of pure dimension let $\HH^i(Y)$ (resp. $\HH^i_c(Y)$) be the $i$-th hypercohomology space 
(resp. hypercohomolgy with compact support) of $Y$ with coefficients in the intersection cohomology complex 
$IC(Y,\bbq)$; let $H^i_c(Y)=H^i_c(Y,\bbq)$. If $K$ is a complex of $l$-adic sheaves on an algebraic variety $Y$ 
we denote by $\ch^iK$ the $i$-th cohomology sheaf of $K$ and by $\ch^i_yK$ the stalk of $\ch^iK$ at $y\in Y$. 
The cardinal of a finite set $S$ is denoted by $|S|$. If $S$ is a set and $f:S@>>>S$ is a map we set 
$S^f=\{s\in S;f(s)=s\}$. We set $\ca=\ZZ[v,v\i]$ where $v$ is an indeterminate.  

\head 1. Preliminaries\endhead
\subhead 1.1\endsubhead
Let $G$ be a connected reductive algebraic group over $\kk$. Let $F':G@>>>G$ be the Frobenius map relative to an 
$F_{q'}$-rational structure on $G$ ($q'$ is a power of $p$). Let $\cb$ the variety of Borel subgroups of $G$. 
Note that $F'$ induces an endomorphism $B\m F'(B)$ of $\cb$. Let $W$ be the Weyl group of $G$ viewed as an 
indexing set for the $G$-orbits on $\cb\T\cb$ (simultaneus conjugation); for $w\in W$ let $\co_w$ be the 
$G$-orbit corresponding to $w$. We regard $W$ as a Coxeter group with set of simple reflections $\{s_i;i\in I\}$ 
in the standard way; let $l:W@>>>\NN$ be the corresponding length function. For any $I'\sub I$ let $s_{I'}$ be 
the longest element of the subgroup of $W$ generated by $\{s_i;i\in I'\}$.
Let $\le$ be the standard partial order on $W$. Now $F'$ induces an automorphism $\d:W@>>>W$ (compatible 
with the length function) by the requirement that $w\in W,(B,B')\in\co_w$ implies $(F'(B),F'(B'))\in\co_{\d(w)}$. Let $\sg:W@>>>\{\pm1\}$ be the homomorphism $w\m(-1)^{l(w)}$.

\subhead 1.2\endsubhead
Let $\Irr W$ be a set of representatives for the isomorphism classes of irreducible representations of $W$ over 
$\QQ$. For $E\in\Irr W$ let $M_E$ be the set of linear maps of finite order $\D:E@>>>E$ such that
$\D(w(e))=\d(w)(\D(e))$ for any $w\in W$, $e\in E$. Then $|M_E|$ is $0$ or $2$. Let
$\Irr_\d W=\{E\in\Irr W;|M_E|=2\}$. 
For $E\in\Irr W$ we define $E^\da\in\Irr W$ by $E^\da\cong E\ot\sg$.

Let $H$ be the Hecke algebra over $\ca(v)$ of $W$ with respect to the weight function $w\m l(w)$ on $W$. Thus $H$
has a $\ca(v)$-basis $(T_w)_{w\in W}$ and we have $T_wT_{w'}=T_{ww'}$ if $w,w'\in W$, $l(ww')=l(w)+l(w')$;
moreover we have $(T_{s_i}-v^2)(T_{s_i}+1)=0$ for $i\in I$. Note that $\QQ[W]=\QQ\ot_\ca H$ where $\QQ$ is 
viewed as an $\ca$-algebra via $v\m1$; $w\in\QQ[W]$ corresponds to $1\ot T_w$. Let $J$ be the ring with 
$\ZZ$-basis $(t_w)_{w\in W}$ defined as in \cite{\LU, 18.3} in terms of $(W,l)$ and let $\Ph:H@>>>\ca\ot J$ be 
the $\ca$-algebra homomorphism defined in \cite{\LU, 18.9}. After applying $\QQ\ot_\ca()$ to $\Ph$ we obtain an 
algebra isomorphism $\Ph_\QQ:\QQ[W]@>\si>>\QQ\ot J$. Via this isomorphism any $E\in\Irr W$ becomes a simple 
$\QQ\ot J$-module denoted by $E_\sp$. Let $E\in\Irr_\d W$, $\D\in M_E$. We define $\D:E_\sp@>>>E_\sp$ by 
$\D(\x)=\Ph_\QQ(\D(\Ph_\QQ\i(\x)))$; we have $t_{\d(w)}(\D(\x))=\D(t_w(\x))$ for any $w\in W,\x\in E_\sp$. After 
applying $\QQ(v)\ot_\ca$ to $\Ph$ we obtain an algebra isomorphism $\Ph_{\QQ(v)}:\QQ(v)\ot H@>\si>>\QQ(v)\ot J$. 
Via this isomorphism the simple $\QQ(v)\ot J$-module $\QQ(v)\ot_\QQ E_\sp$ ($E$ as above) becomes a simple 
$\QQ(v)\ot H$-module denoted by $E_{v^2}$ and the isomorphism $1\ot\D:\QQ(v)\ot_\QQ E_\sp@>>>\QQ(v)\ot_\QQ E_\sp$
becomes an isomorphism $\D: E_{v^2}@>>>E_{v^2}$ such that $T_{\d(w)}(\D(\x))=\D(T_w(\x))$ for any 
$w\in W,\x\in E_{v^2}$. 

We assume that for each $E\in\Irr_\d W$ we have choosen an element $\D\in M_E$.

Let $a:W@>>>\NN$ be the function defined (in terms of $l:W@>>>\NN$) in \cite{\ORA, (5.27.1)} or equivalently as 
in \cite{\LU, 13.6}. Let $E\m a_E$ be the function $\Irr W@>>>\NN$  defined in \cite{\ORA, (4.1.1)}. For 
$E\in\Irr W$ we set $a'_E=a_{E^\da}$.

For $w\in W$ we have
$$v^{-l(w)}\tr(\D T_w,E_{v^2})=c'_{w,\D,E}v^{a'_E}+\text{lower powers of $v$}$$
where $c'_{w,\D,E}\in\ZZ$ (see \cite{\ORA, (5.1.23)}); by an argument similar to that in \cite{\LU, 20.10}, we 
have $c'_{w,\D,E}=\tr(\D',E^\da_\sp)$ where $\D':E^\da@>>>E^\da$ is $\D\ot1:E\ot\sg@>>>E\ot\sg$.

Let $\car_W$ be the vector space of formal linear combinations $\sum_{E\in\Irr_\d W}r_EE$, $r_E\in\QQ$. For 
$w\in W$ we set
$$\fA_w=\sum_{E\in\Irr_\d W}c'_{w,\D,E}E\in\car_W.$$
(Compare \cite{\ORA, (5.11.6)}.)

\subhead 1.3\endsubhead
Let $w\in W$. Following \cite{\DL} we set $X_w=\{B'\in\cb;(B',F'(B'))\in\co_w\}$. Let
$$\bX_w=\{B'\in\cb;(B',F'(B'))\in\cup_{z\in W;z\le w}\co_z\}$$ 
be the closure of $X_w$ in $\cb$. For $x\in G^{F'}$, $Ad(x):\cb@>>>\cb$ (conjugation by $x$) leaves $X_w,\bX_w$ 
stable and induces linear automorphisms $Ad(x)^*$ of $\HH^i(\bX_w),H^i_c(X_w)$. Then $x\m\Ad(x)^*{}\i$ makes 
$\HH^i(\bX_w),H^i_c(X_w)$ into $G^{F'}$-modules. Let $\ce$ be a set of representatives for the isomorphism 
classes of irreducible representations of $G^{F'}$ which appear in $H^i_c(X_w)$ for some $w\in W,i\in\NN$ or 
equivalently in $\HH^i(\bX_w)$ for some $w\in W,i\in\NN$ (thus $\ce$ is the set of unipotent representations of 
$G^{F'}$.) Let $[\ce]_\ZZ$ be the Grothendieck group of the category of representations of $G^{F'}$ which are 
finite sums of unipotent representations. Then $[\ce]:=\QQ\ot[\ce]_\ZZ$ has a basis $\{\r;\r\in\ce\}$. For 
$f\in[\ce]$ let $(\r:f)\in\QQ$ be the coefficient of $\r\in\ce$ in $f$.

For $w\in W$ let $R_w=\sum_i(-1)^iH^i_c(X_w)$ viewed as an element of $[\ce]$. As in \cite{\ORA, 3.7}, for any 
$E\in\Irr_\d W$ we define 
$$R_E=|W|\i\sum_{w\in W}\tr(\D w,E)R_w\in[\ce].$$
For any $\x=\sum_{E\in\Irr_\d W}r_EE\in\car_W$ we set $R_\x=\sum_{E\in\Irr_\d W}r_ER_E\in[\ce]$. In particular, 
for $w\in W$, $R_{\fA_w}$ is defined.

\subhead 1.4\endsubhead
Let $\le_{LR}$ be the preorder on $W$ defined in \cite{\KL} and let $\si_{LR}$ be the corresponding equivalence 
relation on $W$ (the equivalence classes are the two-sided cells of $W$; they form a set $\cc$). For $w,w'\in W$ 
we write $w'<_{LR}w$ when $w'\le_{LR}w$ and $w'\not\si_{LR}w$. The relation $\le_{LR}$ on $W$ induce a relation 
on $\cc$ denoted again by $\le_{LR}$. It is known \cite{\KL} that if $c\in\cc$ then $c^*:=cs_I=s_Ic\in\cc$ and 
that $c\m c^*$ reverses the preorder $\le_{LR}$.

If $E\in\Irr W$ then there is a unique $c\in\cc$ such that $t_w:E_\sp@>>>E_\sp$ is nonzero for some $w\in c$ and 
is zero for all $w\in W-c$; we then write $E\dsv c$. From the definitions we see that

(a) {\it If $c\in\cc,w\in c$ then $\fA_w$ is a linear combination of $E$ such that $E\dsv c^*$.}
\nl
The folowing result is contained in \cite{\ORA, 4.23}.
\proclaim{Theorem 1.5}Let $\r\in\ce$. There exists a unique $c\in\cc$ such that $(\r:R_E)\ne 0$ for some 
$E\in\Irr_\d W$ with $E\dsv c$ and $(\r:R_E)=0$ for all $E\in\Irr_\d W$ with $E\not\dsv c$. We then write 
$\un{\r}=c$.
\endproclaim
Using the theorem and 1.4(a) we see that for $c\in\cc,w\in c$ we have

(a) $R_{\fA_w}\in\sum_{\r\in\ce;\un{\r}=c^*}\QQ\r$.
\nl
The folowing result is contained in \cite{\ORA, 6.15}.
\proclaim{Theorem 1.6} Let $w\in W$. We have
$$\HH^{l(w)-a(w)}(\bX_w)=(-1)^{l(w)-a(w)}R_{\fA_w}+\sum_{w'\in W;w'<_{LR}w}n_{w',w}R_{\fA_{w'}}\in[\ce]\tag a$$
where $n_{w',w}\in\QQ$. 
\endproclaim
For $w\in W,i\in\NN$ and $\r\in\ce$ we denote by $\HH^i(\bX_w)_\r$ the $\r$-isotypic component of the 
$G^{F'}$-module $\HH^i(\bX_w)$ and we set 
$$\HH^{l(w)-a(w)}(\bX_w)_\di=\sum_{\r;\un{\r}=c^*}\HH^{l(w)-a(w)}(\bX_w)_\r$$
where $c\in\cc$ is defined by $w\in c$.

Let $c\in\cc$. Using the theorem and 1.5(a) we see that if $w\in c$ then
$$\HH^{l(w)-a(w)}(\bX_w)\text{ is a $\QQ$-linear combination of $\r\in\ce$ such that $c^*\le_{LR}\un{\r}$}.\tag b
$$
Hence projecting the equality (a) onto the subspace generated by the $\r\in\ce$ such that $\un{\r}=c^*$ we see 
that for any $w\in c$ we have 
$$\HH^{l(w)-a(w)}(\bX_w)_\di=(-1)^{l(w)-a(w)}R_{\fA_w}.\tag c$$
(This projection maps each $R_{\fA_{w'}}$ with $w'<_{LR}w$ to zero, see 1.5(a).) In particular, for any $w\in W$,
$$(-1)^{l(w)-a(w)}R_{\fA_w}\text{ is an actual $G^{F'}$-module}.\tag d$$
\proclaim{Corollary 1.7} Let $\r\in\ce$. We have $\HH^{l(w)-a(w)}(\bX_w))_\r\ne0$ for some $w\in\un{\r}^*$. If 
$c'\in\cc,w'\in c$ and $\HH^{l(w')-a(w')}(\bX_{w'})_\r\ne0$ then $\un{\r}^*\le_{LR}c'$. 
\endproclaim
Let $c=\un{\r}$. 
We can find $E\in\Irr_\d W$ such that $(\r:R_E)\ne0$, $E\dsv c$. By \cite{\ORA, 5.13(ii)}, $E$ is a
$\QQ$-linear combination of elements $\fA_x$ ($x\in c^*$). Hence $R_E$ is a $\QQ$-linear combination of elements 
$R_{\fA_x}$ ($x\in c^*$). It follows that $(\r:\fA_x)\ne0$ for some $x\in c^*$. Hence by 1.6(c) we have 
$\HH^{l(x)-a(x)}(\bX_x))_\r\ne0$. The last sentence in the corollary follows from 1.6(b).

\head 2. The main results\endhead
\subhead 2.1\endsubhead
We preserve the setup of 1.1. In addition, we fix an automorphism $\s:G@>>>G$ such that for some 
$(B,B^*)\in\co_{s_I}$ we have $\s(B)=B$, $\s(B^*)=B^*$ and such that $\s F'=F'\s:G@>>>G$. We also fix an integer 
$m\ge1$ such that $\s^m=1$. Now $\s$ induces a (length preserving) automorphism of $W$ denoted again by $\s$; 
thus, for $w\in W,(B_1,B_2)\in\co_w$ we have $(\s(B_1),\s(B_2))\in\co_{\s(w)}$; moreover $\s\d=\d\s:W@>>>W$. For 
$i\in I$ we have $\s(s_i)=s_{\s(i)}$ where $\s:I@>>>I$ is a bijection. Let 
$\hG$ be the semidirect product of $G$ with the cyclic group of order $m$ with generator $\s$ so that in $\hG$ we
have $\s g\s\i=\s(g)$ for all $g\in G$. Note that $\hG$ is naturally an affine algebraic group with identity
component $G$. We extend $F'$ to a homomorphism $\hG@>>>\hG$ (denoted again by $F'$) by $\s^ig\m\s^iF'(g)$ for 
$i\in[0,m-1],g\in G$. This is the Frobenius map for an $F_{q'}$-rational structure on $\hG$. Note that $\hG^{F'}$
is the semidirect product of $G^{F'}$ with the cyclic group of order $m$ with generator $\s$.

Let $d\ge1$ be an integer with the following property: there exists $(B,B^*)\in\co_{s_I}$ such that $\s(B)=B$, 
$\s(B^*)=B^*$, $F'{}^d(B)=B$, $F'{}^d(B^*)=B^*$, $F'{}^d$ acts on $B\cap B^*$ as $t\m t^{q'{}^d}$ (clearly, such 
$d$ exists). Let $r\ge1$ be a multiple of $d$ and let $F=F'{}^r\s:G@>>>G$, a Frobenius map relative to an 
$F_q$-rational structure on $G$ where $q=q'{}^r$. Note that $w\in W,(B_1,B_2)\in\co_w$ implies 
$(F(B_1),F(B_2))\in\co_{\s(w)}$. 

Let $\dW=\{w\in W;\s(w)=w\}$. It is known that $\dW$ is itself a Weyl group with standard generators $s_\o$ where
$\o$ runs over the set $\dI$ of $\s$-orbits on $I$. Let $\dl:\dW@>>>\NN$ be the length function of $\dW$; thus 
$\dl(s_\o)=1$ for any $\o\in\dI$. The restriction $l|_{\dW}$ of $l$ to $\dW$ is a {\it weight function} (in the 
sense of \cite{\LU, 3.1}) on the Coxeter group $\dW$. We define $\dsg:\dW@>>>\{\pm1\}$ by $\dsg(w)=(-1)^{\dl(w)}$.

\subhead 2.2\endsubhead
Let $\Irr\dW$ be a set of representatives for the isomorphism classes of irreducible representations of $\dW$ 
over $\QQ$. Now $\d:W@>>>W$ restricts to an automorphism of $\dW$ (denoted again by $\d$) preserving the set 
$\{s_\o,\o\in\dI\}$ and the weight function $l|_{\dW}$. 
For $E\in\Irr\dW$ let $M_E$ be the set of linear maps of finite order $\D:E@>>>E$ such that
$\D(w(e))=\d(w)(\D(e))$ for any $w\in\dW$, $e\in E$. Then $|M_E|$ is $0$ or $2$. Let
$\Irr_\d\dW=\{E\in\Irr\dW;|M_E|=2\}$. 
For $E\in\Irr\dW$ we define $E^\da\in\Irr\dW$ by $E^\da\cong E\ot\dsg$.

Let $\dH$ be the Hecke algebra over $\ca$ of $\dW$ with respect to the weight function $l|_{\dW}$. Thus $\dH$ has
an $\ca$-basis $(\dT_w)_{w\in\dW}$ and we have $\dT_w\dT_{w'}=\dT_{ww'}$ if $w,w'\in\dW$, 
$\dl(ww')=\dl(w)+\dl(w')$; 
moreover we have $(\dT_{s_\o}-v^{2l(s_\o)})(\dT_{s_\o}+1)=0$ for $\o\in\dI$. Note that $\QQ[\dW]=\QQ\ot_\ca\dH$ 
where $\QQ$ is viewed as an $\ca$-algebra via $v\m1$; $w\in\QQ[W]$ corresponds to $1\ot\dT_w$. 

Let $\bar{}:\ca@>>>\ca$ be the ring involution such that $v\m v\i$. Let $\bar{}:\dH@>>>\dH$ be the ring 
homomorphism such that $\ov{a\dT_w}=\bar a\dT_{w\i}\i$ for $w\in\dW,a\in\ca$. For any $x,y\in\dW$ we define a 
polynomial $\dR_{x,y}(X)\in\ZZ[X]$, ($X$ is an indeterminate) by the equality
$$\dT_y=\sum_{x\in\dW}\dR_{x,y}(v^2)v^{2l(x)}\dT_{x\i}\i.$$
(See \cite{\KL,(2.0.a)},\cite{\LU, 4.3}.) Note that $\dR_{x,y}=0$ unless $x\le y$. It follows that 
$$\dT_y\dT_{s_I}=\sum_{x\in\dW;x\le y}\dR_{x,y}(v^2)v^{2l(x)}\dT_{xs_I}.\tag a$$
For any $w\in\dW$ there is a unique element $\dC_w\in\dH$ such that
$$\align&\dC_w=\sum_{y\in\dW;y\le w}\dsg(yw)\dP_{y,w}(v^{-2})v^{l(w)-2l(y)}\dT_y\\&=
\sum_{y\in\dW;y\le w}\dsg(yw)\dP_{y,w}(v^2)v^{-l(w)+2l(y)}\dT_{y\i}\i\endalign$$
where $\dP_{y,w}(X)\in\ZZ[X]$ has degree $\le(l(w)-l(y)-1)/2$ if $y<w$ and $\dP_{w,w}(X)=1$. (See 
\cite{\LU,\S5}; compare \cite{\KL}.) For $y\le w$ in $\dW$ we have from the definitions
$$v^{2l(w)}\dP_{y,w}(v^{-2})=\sum_{z\in\dW;y\le z\le w}v^{2l(y)}\dR_{y,z}(v^2)\dP_{z,w}(v^2).\tag b$$
We set $\dP_{y,w}=0$ if $y\not\le w$. We define for $y,w$ in $\dW$ a polynomial $\dQ_{y,w}(X)\in\ZZ[X]$ by the 
requirement that for any $y,w$ in $\dW$,
$$\sum_{z\in\dW}\dsg(zw)\dP_{y,z}(X)\dQ_{z,w}(X)\text{ is $1$ if $y=w$ and is $0$ if $y\ne w$}.\tag c$$
Let $\dJ$ be the ring with $\ZZ$-basis $(\dt_w)_{w\in\dW}$ defined as in \cite{\LU, 18.3} in terms of $\dW$ and 
the weight function $l|_{\dW}$ and let $\Ph:\dH@>>>\ca\ot\dJ$ be the $\ca$-algebra homomorphism defined in 
\cite{\LU, 18.9}. (The definitions and results of \cite{\LU} that were just quoted are applicable in view of 
\cite{\LU, \S16}.) After applying $\QQ\ot_\ca()$ to $\Ph$ we obtain an 
algebra isomorphism $\Ph_\QQ:\QQ[\dW]@>\si>>\QQ\ot\dJ$. Via this isomorphism any $E\in\Irr\dW$ becomes a simple 
$\QQ\ot J$-module denoted by $E_\sp$. Let $E\in\Irr_\d\dW$, $\D\in M_E$. We define $\D:E_\sp@>>>E_\sp$ by 
$\D(\x)=\Ph_\QQ(\D(\Ph_\QQ\i(\x)))$; we have $\dt_{\d(w)}(\D(\x))=\D(\dt_w(\x))$ for any $w\in\dW,\x\in E_\sp$. 
After applying $\QQ(v)\ot_\ca$ to $\Ph$ we obtain an algebra isomorphism
$\Ph_{\QQ(v)}:\QQ(v)\ot_\ca\dH@>\si>>\QQ(v)\ot_\ZZ\dJ$. Via this isomorphism the simple $\QQ(v)\ot\dJ$-module 
$\QQ(v)\ot_\QQ E_\sp$ ($E$ as above) becomes a simple $\QQ(v)\ot_\ca\dH$-module denoted by $E_{v^2}$ and the 
isomorphism $1\ot\D:\QQ(v)\ot_\QQ E_\sp@>>>\QQ(v)\ot_\QQ E_\sp$ becomes an isomorphism $\D: E_{v^2}@>>>E_{v^2}$ 
such that $\dT_{\d(w)}(\D(\x))=\D(\dT_w(\x))$ for any $w\in\dW,\x\in E_{v^2}$. 

We assume that for each $E\in\Irr_\d\dW$ we have choosen an element $\D\in M_E$.
\subhead 2.3\endsubhead
Let $w\in\dW$. Now $X_w$ and $\bX_w$ are stable under $\ph:=F'{}^d:\cb@>>>\cb$ (since $\ph$ commutes with $F'$ 
and it acts trivially on $W$) and under $\s:\cb@>>>\cb$ (since $\s$ commutes with $F'$ and $\s(w)=w$); hence they
are $F$-stable and there are induced automorphisms $\ph^*,(F'{}^r)^*,\s^*,F^*$ of $\HH^i(\bX_w)$ and 
$H^i_c(X_w)$. Hence $X_w,\bX_w$ are stable under the $\hG^{F'}$-action $\s^ig:B\m\s^i(\Ad(g)B)$ 
($i\in[0,m-1],g\in G^{F'}$) and $\hG^{F'}$ acts on $\HH^i(\bX_w)$ and $H^i_c(X_w)$ by $\s^ig\m(\s^*)\i\Ad(g\i)^*$.
Let $\dce$ be a set of representatives for the isomorphism classes of irreducible representations of $\hG^{F'}$
whose restriction to $G^{F'}$ is a direct sum of unipotent representations. Let $\ce_0$ be the set of all 
$\r\in\ce$ such that $\r$ extends to a $\hG^{F'}$-module; let $\dce_0$ be the set of all $\dr\in\dce$ such that 
$\dr|_{G^{F'}}$ is irreducible. Let $[\dce]\ti{}_\ZZ$ be the Grothendieck group of the category of 
representations of $\hG^{F'}$ whose restriction to $G^{F'}$ are finite sums of unipotent representations. Then
$[\dce]\ti{}:=\QQ\ot[\dce]\ti{}_\ZZ$ has a basis $\{\dr;\dr\in\dce\}$. For $\dr\in\dce$, $f\in[\dce]\ti{}$ we 
denote by $(\dr:f)\ti{}$ the coefficient of $\dr$ in $f$. Let $[\dce]$ be the quotient of $[\dce]\ti{}$ by the 
subspace consisting of all elements whose character restricted to $G^{F'}\s$ is zero. Note that for $h\in[\dce]$,
$\tr(x\s,h)$ makes sense for any $x\in G^{F'}$ and we can define for $\dr\in\dce$, 
$$(\dr:h)=|G^{F'}|\i\sum_{x\in G^{F'}}\tr((x\s)\i,\dr)\tr(x\s,h).$$
We show that if $f\in[\dce]\ti{}$, $h$ is its image in $[\dce]$ and $\dr\in\dce$, then
$$(\dr:h)=\sum_{\c\in L}(\dr\ot\c:f)\ti{}\c(\s)\i\tag a$$
where $L$ is the set of one dimensional characters of $\hG^{F'}$, trivial on $G^{F'}$. Indeed,
$$\align&\sum_{\c\in L}(\dr\ot\c:f)\ti{}\c(\s)\i\\&
=|G^{F'}|\i m\i\sum_{\c\in L}\sum_{x\in G^{F'},j\in[0,m-1]}\tr(x\s^j,\dr\ot\c)\tr((x\s^j)\i,f)\c(\s)\i\\&=
|G^{F'}|\i m\i\sum_{x\in G^{F'},j\in[0,m-1]}\tr(x\s^j,\dr)\tr((x\s^j)\i,f)\sum_{\c\in L}\c(\s^{j-1})\\&
=|G^{F'}|\i\sum_{x\in G^{F'}}\tr(x\s,\dr)\tr((x\s)\i,f)=(\dr:h),\endalign$$
as desired.

For any $\r\in\ce_0$ let $[\dce]^\r$ be the (one dimensional) subspace of $[\dce]$ generated by the elements 
$\dr\in\dce_0$ such that $\r=\dr|_{G^{F'}}$. From the definitions we have a direct sum decomposition
$$[\dce]=\op_{\r\in\ce_0}[\dce]^\r.$$
For any $c\in\cc$ we set $[\dce]^c=\op_{\r\in\ce_0;\un{\r}^*=c}[\dce]^\r$. We have 
$$[\dce]=\op_{c\in\cc}[\dce]^c.$$
The following result is clear from the definitions.

(b) {\it Let $h\in[\dce]$ be such that $(\dr:h)=0$ for any $\dr\in\dce$. Then $h=0$.}
\nl
For $w\in\dW$ let 
$$\dR_w=\sum_i(-1)^iH^i_c(X_w)$$ 
viewed as an element of $[\dce]$. For any $E\in\Irr_\d\dW$ we set
$$\dR_E=|\dW|\i\sum_{w\in\dW}\tr(\D w,E)\dR_w\in[\dce].$$
For any $\x\in\dH$ any $E\in\Irr_\d\dW$ and any $i\in\ZZ$ we define $\tr(\D h,E_{v^2};i)\in\ZZ$ by
$$\tr(\D h,E_{v^2})=\sum_i\tr(\D h,E_{v^2};i)v^i.$$
Part (i) of the following theorem is a partial generalization of \cite{\ORA, 3.8(ii)}; part (ii) (a partial
generalization of \cite{\ORA, 4.23}, see also 1.5) provides an affirmative answer to a conjecture of G. Malle 
\cite{\MA} (which he proved in the case where $\rank(G)\le 6$).
\proclaim{Theorem 2.4} (i) For any $w\in\dW$ the following equality holds in $\QQ[v]\ot[\dce]$: 
$$\align&\sum_i(-1)^i\HH^i(\bX_w)v^i=\sum_{E\in\Irr_\d\dW}\tr(\D\sum_{y\in\dW}\dP_{y,w}(v^2)\dT_y,E_{v^2})\dR_E.
\tag a\endalign$$
Equivalently, for any $i\in\ZZ$ we have
$$(-1)^i\HH^i(\bX_w)=\sum_{E\in\Irr_\d\dW}\tr(\D\sum_{y\in\dW}\dP_{y,w}(v^2)\dT_y,E_{v^2};i)\dR_E.$$
(ii) Let $E\in\Irr_\d\dW$. There exists $c\in\cc$ such that $\dR_E\in[\dce]^c$.
\endproclaim
The proof of (i) (which has much in common with that of \cite{\ORA, 3.8(ii)}) is given in 2.15. The proof of (ii) 
is given in 2.19 where the two-sided cell $c$ in (ii) is explicitly described in terms of $E$ (see 2.19(c)).

\subhead 2.5\endsubhead
We fix a square root $q'{}^{1/2}$ of $q'$ in $\bbq$; for any integer $n$ we set $q'{}^{n/2}=(q'{}^{1/2})^n$,
$q^{n/2}=q'{}^{nr/2}$. Let $\dH_q=\bbq\ot_\ca\dH$ where $\bbq$ is viewed as an $\ca$-algebra via $v\m q^{1/2}$.

Let $\cf$ be the vector space of functions $\cb^F@>>>\bbq$. For any $w\in\dW$ we define a linear map
$\dT_w:\cf@>>>\cf$ by $\dT_w(f)=f'$ for $f\in\cf$ where $f'(B)=\sum_{B'\in\cb^F;(B,B')\in\co_w}f(B')$ 
for $B\in\cb^F$. In this way $\cf$ becomes a $\dH_q$-module. Note that the linear maps $\dT_w:\cf@>>>\cf$ 
($w\in\dW$) are linearly independent. For $w,w'\in\dW$ we have
$\dT_w\dT_{w'}=\sum_{w''\in\dW}N_{w,w',w''}\dT_{w''}$ where
$$N_{w,w',w''}=|\{B\in\cb^F;(B_1,B)\in\co_w,(B,B_2)\in\co_{w'}\}|$$
for any $(B_1,B_2)\in\co_{w''}^F$. Using this and 2.2(a) we see that
$$\sum_{z\in\dW}N_{y,s_I,z}\dT_z=\sum_{x\in\dW}\dR_{x,y}(q)q^{l(x)}\dT_{xs_I}$$
for any $y\in\dW$ as linear maps $\cf@>>>\cf$. Using the linear independence of the linear maps $\dT_z$ we deduce
$$\dR_{x,y}(q)q^{l(x)}=N_{y,s_I,xs_I}\tag a$$
for any $x,y$ in $\dW$.

\subhead 2.6\endsubhead
Now let $B,B^*$ be as in 2.1. For any $w\in\dW$ we define $B_w\in\cb$ by the conditions
$(B,B_w)\in\co_w,(B_w,B^*)\in\co_{w\i s_I}$. Note that $B_w$ is fixed by $F'{}^r:\cb@>>>\cb$ and by 
$\s:\cb@>>>\cb$ hence by $F:\cb@>>>\cb$. For $z\in\dW$ let 
$$\cb_z=\{B'\in\cb;(B,B')\in\co_z\},$$ 
$$\bcb_z=\{B'\in\cb;(B,B')\in\cup_{z'\in W;z'\le z}\co_z\},$$
$$A^z=\{B'\in\cb;(B',B_{zs_I})\in\co_{s_I}\};$$
note that $\cb_z,\bcb_z,A^z$ are stable under $F'{}^r:\cb@>>>\cb$ and under $\s:\cb@>>>\cb$ hence under
$F:\cb@>>>\cb$. 

Let $w\in\dW$. Let $K_w=IC(\bcb_w,\bbq)$ be the intersection cohomology complex of $\bcb_w$. For any $y\in\dW$ 
such that $y\le w$, $\s$ induces an isomorphism of local systems 
$\s^*(\ch^iK_w|_{\cb_y})@>\si>>\ch^iK_w|_{\cb_y}$; this induces an automorphism 
$\s^*:\ch^i_{B_y}K_w@>\si>>\ch^i_{B_y}K_w$ since $B_y\in\cb_y$ is fixed by $\s$. The trace of this automorphism 
is denoted by $n_{i,y,w,\s}$. 

The following generalization of \cite{\KLL, 4.3} was stated without proof in \cite{\LC, (8.1)}; its proof (which 
uses \cite{\KLL, 4.2}) is a generalization of that of \cite{\KLL, 4.3}.
\proclaim{Theorem 2.7} For any $y,w\in\dW$ such that $y\le w$ we have
$$\dP_{y,w}(X)=\sum_{i\in\NN}n_{2i,y,w,\s}X^i.$$
\endproclaim
From the definitions, for $z,y\in\dW$ we have $|(\cb_z\cap A^y)^F|=N_{z,s_I,ys_I}$; using 2.5(a), it follows that
$$|(\cb_z\cap A^y)^F|=\dR_{y,z}(q)q^{l(y)}.\tag a$$
Note that $K_w|_{\bcb_w\cap A^y}$ is the intersection cohomology complex $IC(\bcb_w\cap A^y,\bbq)$ of 
$\bcb_w\cap A^y$ since $\bcb_w\cap A^y$ is open in $\bcb_w$. Now the Grothendieck-Lefschetz trace formula for 
$F:\bcb_w\cap A^y@>>>\bcb_w\cap A^y$ gives:
$$\sum_i(-1)^i\tr(F^*,\HH^i_c(\bcb_w\cap A^y))=\sum_{B'\in(\bcb_w\cap A^y)^F}\sum_i(-1)^i\tr(F^*,\ch^i_{B'}K_w).$$
By specifying $z\in W$ such that $B'\in\co_z^F$ (so that $z$ is necessarily in $\dW$) we obtain
$$\align&\sum_i(-1)^i\tr(F^*,\HH^i_c(\bcb_w\cap A^y))\\&=\sum_{z\in\dW;z\le w}\sum_{B'\in(\cb_z\cap A^y)^F}
\sum_i(-1)^i\tr(F^*,\ch^i_{B'}K_w).\endalign$$  
Note that $\tr(F^*,\ch^i_{B'}K_w)$ is independent of $B'$ when $B'$ runs through $(\cb_z\cap A^y)^F$ and even 
when $B'$ runs through $\cb_z^F$ (since the local system $H^iK_w|_{\cb_z}$ is $B$-equivariant for the obvious 
transitive $B$-action on $\cb_z$ with connected isotropy groups); moreover we have $B_z\in\cb_z^F$ and we deduce 
that
$$\align&\sum_i(-1)^i\tr(F^*,\HH^i_c(\bcb_w\cap A^y))\\&=\sum_{z\in\dW;z\le w}|(\cb_z\cap A^y)^F|
\sum_i(-1)^i\tr(F^*,\ch^i_{B_z}K_w);\endalign$$
using (a) and the fact that $\dR_{y,z}(X)=0$ unless $y\le z$, we deduce that
$$\align&\sum_i(-1)^i\tr(F^*,\HH^i_c(\bcb_w\cap A^y))\&=\sum_{z\in\dW;y\le z\le w}\dR_{y,z}(q)q^{l(y)}
\sum_i(-1)^i\tr(F^*,\ch^i_{B_z}K_w).\endalign$$
By Poincar\'e duality for intersection cohomology we have
$$\sum_i(-1)^i\tr(F^*,\HH^i_c(\bcb_w\cap A^y))=q^{l(w)}\sum_i(-1)^i\tr(F^*{}\i,\HH^i(\bcb_w\cap A^y))$$
since $\bcb_w\cap A^y$ is of pure dimension $l(w)$ (or is empty). By \cite{\KLL, 1.5, 4.5} we have
$$\tr(F^*{}\i,\HH^i(\bcb_w\cap A^y))=\tr(F^*{}\i,\ch^i_{B_y}K_w)$$
hence
$$\align&q^{l(w)}\sum_i(-1)^i\tr(F^*{}\i,\ch^i_{B_y}K_w)\\&=\sum_{z\in\dW;y\le z\le w}\dR_{y,z}(q)q^{l(y)}
\sum_i(-1)^i\tr(F^*,\ch^i_{B_z}K_w).\endalign$$
By \cite{\KLL, 4.2}, $\ch^i_{B_z}K_w$ is zero if $i$ is odd while if $i$ is even, $(F'{}^r)^*$ acts on 
$\ch^i_{B_z}K_w$ as $q^{i/2}$ times a unipotent transformation hence $F^*$ acts on $\ch^i_{B_z}K_w$ as $q^{i/2}$ 
times a unipotent transformation times the action of $\s^*$ (which commutes with the unipotent transformation, 
since $\s F'{}^r=F'{}^r\s$). Thus we have
$$\align&q^{l(w)}\sum_{i\in2\NN}q^{-i/2}\tr(\s^*{}\i,\ch^i_{B_y}K_w)\\&
=\sum_{z\in\dW;y\le z\le w}\dR_{y,z}(q)q^{l(y)}\sum_{i\in2\NN}q^{i/2}\tr(\s^*,\ch^i_{B_z}K_w).\endalign$$
Now for any $z\in\dW$ such that $z\le w$ we set:
$$\tP_{z,w}(X)=\sum_{i\in2\NN}\tr(\s^*,\ch^i_{B_z}K_w)X^{i/2}\in R[X],$$
$$\tP'_{z,w}(X)=\sum_{i\in2\NN}\tr(\s^*{}\i,\ch^i_{B_z}K_w)X^{i/2}\in R[X]$$
where $R$ is the subring of $\bbq$ generated by the $m$-th roots of $1$. By the definition of intersection 
cohomology, $\tP_{z,w},\tP'_{z,w}$ are polynomials in $X$ of degree $\le(l(w)-l(z)-1)/2$ (if $z<w$) and are equal
to $1$ if $z=w$. We have
$$q^{l(w)}\tP'_{y,w}(q\i)=\sum_{z\in\dW;y\le z\le w}q^{l(y)}\dR_{y,z}(q)\tP_{z,w}(q).$$
Since here $q$ is an arbitrary power of $q'{}^d$ it follows that
$$X^{l(w)}\tP'_{y,w}(X\i)=\sum_{z\in\dW;y\le z\le w}X^{l(y)}\dR_{y,z}(X)\tP_{z,w}(X).\tag b$$
Note also that
$$X^{l(w)}\dP_{y,w}(X\i)=\sum_{z\in\dW;y\le z\le w}X^{l(y)}\dR_{y,z}(X)\dP_{z,w}(X).\tag c$$
We show by induction on $l(w)-l(y)$ that
$$\tP_{y,w}(X)=\tP'_{y,w}(X)=\dP_{y,w}(X)\tag d$$
for any $y\in\dW$ such that $y\le w$. If $l(y)=l(w)$ we have $y=w$ and all three terms in (d) are $1$. Now assume
that $y<w$ and that the result is known when $y$ is replaced by $z\in\dW$ where $y<z\le w$. Substracting (c) from
(b) we then find (using that $\dR_{y,y}=1$):
$$X^{l(w)}(\tP'_{y,w}(X\i)-\dP_{y,w}(X\i))=X^{l(y)}(\tP_{y,w}(X)-\dP_{y,w}(X)).$$
Thus
$$X^{(l(w)-l(y)/2}(\tP'_{y,w}(X\i)-\dP_{y,w}(X\i))=X^{(-l(w)+l(y))/2}(\tP_{y,w}(X)-\dP_{y,w}(X)).\tag e$$
The left hand side of (e) belongs to $\sum_{n>0}RX^{-n/2}$; the right hand side of (e) belongs to 
$\sum_{n>0}RX^{n/2}$. Hence both sides must be zero. This proves (d). The theorem is proved.

\subhead 2.8\endsubhead
For $w\in\dW$ let $L_w=IC(\bX_w,\bbq)$ be the intersection complex of $\bX_w$ (a variety of pure dimension 
$l(w)$). For any $x\in G^{F'}$ and any $B'\in\bX_w$ such that $Ad(x)F(B')=B'$ we have an induced isomorphism
$F^*\Ad(x)^*:\ch^i_{B'}L_w@>>>\ch^i_{B'}L_w$. We show:

\proclaim{Proposition 2.9} In the setup of 2.8 we have 
$$\align&\sum_i(-1)^i\tr(F^*\Ad(x)^*,\HH^i(\bX_w))\\&=\sum_{y\in\dW;y\le w}\dP_{y,w}(q)\sum_i(-1)^i
\tr(F^*\Ad(x)^*,H^i_c(X_y)).\endalign$$
\endproclaim
By the Grothendieck-Lefschetz trace formula we have
$$\align&\sum_i(-1)^i\tr(F^*\Ad(x)^*,\HH^i(\bX_w))\\&
=\sum_{B'\in\bX_w;Ad(x)F(B')=B'}\sum_i(-1)^i\tr(F^*\Ad(x)^*,\ch^i_{B'}L_w).\tag a\endalign$$
(Note that $Ad(x)F=\Ad(x)\s F'{}^r:\bX_w@>>>\bX_w$ is a Frobenius map relative to an $F_q$-rational structure 
since $Ad(x)\s$ is an automorphism of finite order of $\bX_w$ commuting with the Frobenius map 
$F'{}^r:\bX_w@>>>\bX_w$.) In the sum over $B'$ in (a) we can specify $y\in W$ such that $B'\in X_y$; we have 
necessarily $y\in\dW$. (Indeed assume that $B'\in X_y$. Then $Ad(x)F(B')\in X_{\s(y)}$; hence if $Ad(x)F(B')=B'$ 
then $\s(y)=y$.) Thus the right hand side of (a) becomes
$$\sum_{y\in\dW;y\le w}\sum_{B'\in X_y;Ad(x)F(B')=B'}\sum_i(-1)^i\tr(F^*\Ad(x)^*,\ch^i_{B'}L_w)$$
and it is enough to show that for any $y\in\dW$ such that $y\le w$ we have
$$\align&\sum_{B'\in X_y;Ad(x)F(B')=B'}\sum_i(-1)^i\tr(F^*\Ad(x)^*,\ch^i_{B'}L_w)\\&=
\dP_{y,w}(q)\sum_i(-1)^i\tr(F^*\Ad(x)^*,H^i(X_y)).\endalign$$
By the Grothendieck-Lefschetz trace formula we have 
$$\sum_i(-1)^i\tr(F^*\Ad(x)^*,H^i(X_y))=|\{B'\in X_y;Ad(x)F(B')=B'\}|.\tag b$$
Thus it is enough to show that for any $B'\in X_y$ such that $Ad(x)F(B')=B'$ we have
$$\sum_i(-1)^i\tr(F^*\Ad(x)^*,\ch^i_{B'}L_w)=\dP_{y,w}(q).$$
Let $\dw\in G$ be such that $\dw(B\cap B^*)\dw\i=B\cap B^*$ and $(B,\dw B\dw\i)\in\co_w$. Let $G_w=B\dw B$, 
$G'_w=\{g\in G; g\i F'(g)\in G_w\}$. Let $\bG_w=\cup_{z\in W;z\le w}G_z$, $\bG'_w=\cup_{z\in W;z\le w}G'_z$ be 
the closures of $G_w,G'_w$ in $G$. Let $\tK_w$ (resp. $\tL_w$) be the intersection cohomology complex of $\bG_w$ 
(resp. $\bG'_w$) with coefficients in $\bbq$. Define $\a:\bG'_w@>>>\bX_w$ by $g\m gBg\i$ (a principal 
$B$-bundle). Define $\Ph:\bG'_w@>>>\bG'_w$ by $\Ph(g)=xF(g)$; note that $\Ph$ is a Frobenius map for an 
$F_q$-rational structure on $\bG'_w$ and $\a(\Ph(g))=Ad(x)F(\a(g))$ for any $g\in\bG'_w$. Hence $\a\i(B')$ is 
$\Ph$-stable. Since $\a\i(B')$ is a homogeneous $B$-space with a compatible $F_q$-rational structure given by 
$\Ph$ we can find $g'\in\a\i(B')$ such that $\Ph(g')=g'$. We have canonically $\ch^i_{B'}L_w=\ch^i_{g'}\tL_w$ and
$\tr(F^*\Ad(x)^*,\ch^i_{B'}L_w)=\tr(\Ph^*,\ch^i_{g'}\tL_w)$; moreover we have $g'\in(G'_y)^\Ph$. Define 
$\cl:\bG'_w@>>>\bG_w$ by $g\m g\i F'(g)$ (an \'etale covering). Define $\Ph':\bG_w@>>>\bG_w$ by $g\m F(g)$; note 
that $\cl(\Ph(g))=\Ph'(\cl(g))$ for any $g\in\bG'_w$. We have canonically $\ch^i_{g'}\tL_w=\ch^i_{\cl(g')}\tK_w$ 
and $\tr(\Ph^*,\ch^i_{g'}\tL_w)=\tr(\Ph'{}^*,\ch^i_{\cl(g)}\tK_w)$; moreover, $\cl(g)\in G_y^{\Ph'}$. Define 
$\b:\bG_w@>>>\bcb_w$ by $g\m gBg\i$ (a principal $B$-bundle). Note that $\b(\Ph'(g))=F(\b(g))$ for any 
$g\in\bG_w$. We have canonically $\ch^i_{\cl(g')}\tK_w=\ch^i_{\b(\cl(g'))}K_w$ and 
$\tr(\Ph'{}^*,\ch^i_{\cl(g)}\tK_w)=\tr(F^*,\ch^i_{\b(\cl(g))}K_w)$; moreover, $\b(\cl(g))\in\cb_y^F$. By the 
proof of 2.7 we have 
$$\tr(F^*,\ch^i_{\b(\cl(g))}K_w)=\tr(F^*,\ch^i_{B_y}K_w)=n_{i,y,w,\s}q^{i/2}.$$
It remains to use the equality
$$\dP_{y,w}(q)=\sum_{i\in\NN}(-1)^in_{i,y,w,\s}q^{i/2}$$
which folows from 2.7. The proposition is proved.

\subhead 2.10\endsubhead
Let $w\in\dW$. The equality in Proposition 2.9 can be rewitten as
$$\align&\sum_i(-1)^i\tr((\ph^*)^{r'}\s^*\Ad(x)^*,\HH^i(\bX_w))\\&=\sum_{y\in\dW;y\le w}\dP_{y,w}(q_0^{r'})
\sum_i(-1)^i\tr((\ph^*)^{r'}\s^*\Ad(x)^*,H^i_c(X_y))\endalign$$
for $r'\in\{1,2,3,\do\}$ where $q_0=(q')^d$. Here we can make formally $r'$ tend to zero and we obtain
$$\align&\sum_i(-1)^i\tr(\s^*\Ad(x)^*,\HH^i(\bX_w))\\&=\sum_{y\in\dW;y\le w}\dP_{y,w}(1)\sum_i(-1)^i
\tr(\s^*\Ad(x)^*,H^i_c(X_y)).\tag a\endalign$$ 

\subhead 2.11\endsubhead
If $x\in G^F$ we have $F'(x)\in G^F$ since $FF'=F'F$; moreover we have $F'{}^{rm}(x)=x$ since $F'{}^{rm}=F^m$. 
Let $\wt{G^F}$ be the semidirect product of $G^F$ with the cyclic group of order $rm$ with generator $\th$ in 
which we have $\th z\th\i=F'(z)$ for all $z\in G^F$. 

For any $z\in G^F$ we define a linear map $\uz:\cf@>>>\cf$ ($\cf$ as in 2.2) by $\uz(f)=f'$ ($f\in\cf$) where 
$f'(B_1)=f(z\i B_1z)$ for $B_1\in\cb^F$. Now $z\m\uz$ makes $\cf$ into a $G^F$-module. We define a linear map 
$\th:\cf@>>>\cf$ by $f\m f'$ ($f\in\cf$) where $f'(B_1)=f(F'{}\i(B_1))$ for $B_1\in\cb^F$ (this is well defined 
since $FF'=F'F$). This linear map together with the $G^F$-module structure define a $\wt{G^F}$-module structure
on $\cf$.

Following Shintani, for each $x\in G^{F'}$ we define (noncanonically) an element $\hx\in G^F$ as follows. We 
write $x=aF(a\i)$ with $a\in G$ and we set $\hx=a\i F'(a)$. Note that $\hx\in G^F$ (since $FF'=F'F$) hence 
$\uhx:\cf@>>>\cf$ is well defined. We have the following result (compare \cite{\ORA, 2.10}):

\proclaim{Proposition 2.12} For any $w\in\dW$ and any $x\in G^{F'}$ we have (with notation of 2.11):
$$\sum_i(-1)^i\tr(F^*x^*,H^i_c(X_w))=\tr(\th\i\uhx\i\dT_{w\i},\cf).$$
\endproclaim
By 2.9(b) with $y$ replaced by $w$, the left hand side of the equality above is equal to
$|\{B'\in X_w;Ad(x)F(B')=B'\}|$ hence to 
$$|\{B'\in\cb;(B',F'(B'))\in\co_w,aF(a\i)F(B')F(a)a\i=B'\}|.$$
Setting $B_1=a\i B'a$, we see that the last number is equal to
$$\align&|\{B_1\in\cb;(aB'_1a\i,F'(a)F'(B_1)F'(a\i))\in\co_w;F(B_1)=B_1\}|\\&=
|\{B_1\in\cb^F;(B_1,\hx F'(B_1)\hx\i\in\co_w\}|.\endalign$$
For any $B_1\in\cb^F$ we denote by $f_{B_1}$ the function $\cb^F@>>>\bbq$ such that $f_{B_1}(B')$ is $1$ if
$B'=B_1$ and $0$ if $B'\ne B_1$. Let $f_1=\dT_{w\i}(f_{B_1})$. We have $f_1(B'')=1$ if $(B'',B_1)\in\co_{w\i}$, 
$f_1(B'')=0$ if $(B'',B_1)\n\co_{w\i}$. Let $f_2=\uhx\i(f_1)$. We have $f_2(B'')=1$ if 
$(\hx B''\hx\i,B_1)\in\co_{w\i}$, $f_2(B'')=0$ if $(\hx B''\hx\i,B_1)\n\co_{w\i}$. Let $f_3=\th\i(f_2)$. We have 
$f_3(B'')=1$ if $(\hx F'(B'')\hx\i,B_1)\in\co_{w\i}$, $f_2(B'')=0$ if $(\hx F'\hx\i(B''),B_1)\n\co_{w\i}$. We see
that 
$$\align&\tr(\th\i\uhx\i\dT_{w\i},\cf)=|\{B_1\in\cb^F;(\hx F'(B_1)\hx\i,B_1)\in\co_{w\i}\}|\\&=
|\{B_1\in\cb^F;(B_1,\hx F'(B_1)\hx\i)\in\co_w\}|.\endalign$$
This completes the proof.

\subhead 2.13\endsubhead
Let $w\in\dw$, $x\in G^{F'}$. Combining 2.9 and 2.12 we obtain
$$\sum_i(-1)^i\tr(F^*\Ad(x)^*,\HH^i(\bX_w))
=\sum_{y\in\dW;y\le w}\dP_{y,w}(q)\tr(\th\i\uhx\i\dT_{y\i},\cf).\tag a$$
By \cite{\ORA, 2.20, 3.8}, for $\r\in\ce$, $i\in\NN$, $\ph^*:\HH^i(\bX_w)@>>>\HH^i(\bX_w)$ (notation of 2.10) 
leaves stable $\HH^i(\bX_w)_\r$ (notation of 1.6) and acts on it as $\l_\r q'{}^{id/2}U$ where $U$ is a unipotent
transformation of $\HH^i(\bX_w)_\r$ and $\l_\r$ is a root of $1$ depending only on $\r$ (not on $w,i$); also if 
$\r\in\ce-\ce_0$, $\s^*$ maps $\HH^i(\bX_w)_\r$ to $\HH^i(\bX_w)_{\r'}$ where $\r'\ne\r$, while if $\r\in\ce_0$ 
then $\s^*$ leaves stable $\HH^i(\bX_w)_\r$ and the $\s^*{}\i$ action makes $\HH^i(\bX_w)_\r$ into a 
$\hG^{F'}$-module. It follows that
$$\align&\sum_i(-1)^i\tr(F^*x^*,\HH^i(\bX_w))\\&
=\sum_i(-1)^i\sum_{\r\in\ce_0}\l_\r^{r/d}q^{i/2}\tr(\s^*x^*,\HH^i(\bX_w)_\r).\tag b\endalign$$
We can find $d'\ge1$ such that $\l_\r^{d'}=1$ for all $\r\in\ce$. From now on we assume that $r$ is a multiple of
$dd'$. Then (b) becomes
$$\align&\sum_i(-1)^i\tr(F^*x^*,\HH^i(\bX_w))\\&=\sum_i(-1)^i\sum_{\dr\in\dce_0}q^{i/2}(\dr:\HH^i(\bX_w))\ti{}
\tr((x\s)\i,\dr).\tag c\endalign$$

\subhead 2.14\endsubhead
After applying $\bbq\ot_\ca()$ to $\Ph$ in 2.2 (where $\bbq$ is viewed as an $\ca$-algebra via $v\m q^{1/2}$) we 
obtain an algebra isomorphism $\Ph_q:\dH_q@>\si>>\bbq\ot\dJ$. If $E\in\Irr\dW$ then the simple 
$\bbq\ot\dJ$-module $\bbq\ot E_\sp$ (see 2.2) can be viewed via $\Ph_q$ as a simple $\dH_q$-module denoted by 
$E_q$. If $E\in\Irr_\d\dW$ we define $\D:E_q@>>>E_q$ to correspond under $\Ph_q$ to 
$1\ot\D:\bbq\ot E_\sp@>>>\bbq\ot E_\sp$; we have $\dT_{\d(w)}(\D(\x))=\D(\dT_w(\x))$ for any 
$w\in\dW,\x\in E_q$.

For $E\in\Irr\dW$ let $\cf_E$ be the $E_q$-isotypic component of the $H_q$-module $\cf$; thus we have
$\cf=\op_{E\in\Irr\dW}\cf_E$. Note that each $\cf_E$ is a $G^F$-submodule of $\cf$. From the definitions, for any
$w\in\dW$ we have $\dT_w\th\i=\th\i\dT_{\d(w)}:\cf@>>>\cf$. It follows that $\th\i:\cf@>>>\cf$ permutes the 
summands $\cf_E$ of $\cf$ according to the permutation of $\Irr\dW$ induced by $\d:W@>>>W$. In particular only 
the summands corresponding to $E\in\Irr_\d\dW$ are mapped into themselves. We see that for $z\in G^F$, $y\in\dW$,
we have
$$\tr(\th\i\uz\i\dT_{y\i},\cf)=\sum_{E\in\Irr_\d\dW}\tr(\th\i\uz\i\dT_{y\i},\cf_E).$$
For $E\in\Irr\dW$ we set $V_E=\Hom_{H_q}(E_q,\cf)$; this is an irreducible $G^F$-module (the $G^F$-module
structure comes from that of $\cf$); we have a canonical isomorphism of $(H_q,G^F)$-modules 
$E_q\ot V_E@>\si>>\cf_E$. Via this isomorphism, assuming $E\in\Irr_\d\dW$, the map $\th\i:\cf_E@>>>\cf_E$ becomes 
an isomorphism $X_E:E_q\ot V_E@>>>E_q\ot V_E$ necessarily of the form $X_E=X'_E\ot X''_E$ where 
$X'_E:E_q@>>>E_q$, $X''_E:V_E@>>>V_E$ are isomorphisms such that $\dT_wX'_E=X'_E\dT_{\d(w)}:E_q@>>>E_q$ 
and $zX''_E=X''_EF'(z):V_E@>>>V_E$ for any 
$z\in G^F$. Thus $X'_E$ acts on $E_q$ as a nonzero constant times $\D\i:E_q@>>>E_q$; we may assume that the 
constant is $1$ (by absorbing it into $X''_E$) so that $X'_E=\D\i$. Since $X_E^{rm}=1$ and $(X'_E)^m=1$ we see 
that $(X''_E)^{rm}=1$ so that the $G^F$-module structure on $V_E$ extends to a $\wt{G^F}$-module structure in 
which $\th\i$ acts as $X''_E$. We see that for $z\in G^F$, $y\in\dW$, we have
$$\tr(\th\i\uz\i\dT_{y\i},\cf)=\sum_{E\in\Irr_\d\dW}\tr(\th\i z\i,V_E)\tr(\D\i\dT_{y\i},E_q).$$
We now substitute this (with $z=\hx$) into 2.13(a), taking into account 2.13(c) and using that
$\tr(\D\i\dT_{y\i},E_q)=\tr(\dT_y\D,E_q)=\tr(\D\dT_y,E_q)$ we obtain
$$\align&\sum_i(-1)^i\sum_{\dr\in\dce_0}q^{i/2}(\dr:\HH^i(\bX_w))\ti{}\tr((x\s)\i,\dr)\\&
=\sum_{y\in\dW}\dP_{y,w}(q)\sum_{E\in\Irr_\d\dW}\tr(\th\i\hx\i,V_E)\tr(\D\dT_y,E_q)\tag a\endalign$$
for any $w\in\dW$, $x\in G^{F'}$.

\subhead 2.15\endsubhead
Multiplying both sides of 2.14(a) by $\dsg(wu)\dQ_{w,u}(q)$ with $u\in\dW$ and summing over all $w\in\dW$ we 
obtain
$$\align&\sum_i(-1)^i\sum_{w\in\dW}\sum_{\dr\in\dce_0}q^{i/2}(\dr:\HH^i(\bX_w))\ti{}
\T\dsg(wu)\dQ_{w,u}(q)\tr((x\s)\i,\dr)\\&=\sum_{E\in\Irr_\d\dW}\tr(\th\i\hx\i,V_E)\tr(\D\dT_u,E_q)\endalign$$
for any $u\in\dW$, $x\in G^{F'}$. Multiplying both sides of the last equality by \lb
$|G^{F'}|\i\tr(x\s,\dr')$ with $\dr'\in\dce_0$ and summing over all $x\in G^{F'}$ we obtain
$$\align&\sum_i(-1)^i\sum_{w\in\dW}\sum_{\c\in L}q^{i/2}(\dr'\ot\c:\HH^i(\bX_w))\ti{}\dsg(wu)\dQ_{w,u}(q)\c(\s)\i
\\&=|G^{F'}|\i\sum_{x\in G^{F'}}\sum_{E\in\Irr_\d\dW}\tr(\th\i\hx\i,V_E)\tr(\D\dT_u,E_q)\tr(x\s,\dr')
\tag a\endalign$$
for any $u\in\dW$, $\dr'\in\dce_0$. Here $L$ is as in 2.3(a). We have used that, if $\dr,\dr'\in\dce_0$, then 
$|G^{F'}|\i\sum_{x\in G^{F'}}\tr((x\s)\i,\dr)\tr(x\s,\dr')$ is $\c(\s)\i$ if $\dr=\dr'\ot\c$ for some 
$\c\in L$ and is $0$, otherwise. Next we note that for $E,E'\in\Irr_\d\dW$, 
$$\sum_{u\in\dW}q^{l(u)}\tr(\D\dT_u,E'_q)\tr(\D\dT_u,E_q)$$ 
is equal to $0$ if $E\ne E'$ and is equal to $D(E)(q)$ for some polynomial $D(E)(X)\in\QQ[X]$ if $E=E'$;
moreover, $D(E)(q)$ is a nonzero rational number.

Multiplying both sides of (a) by 
$D(E')(q)\i q^{l(u)}\tr(\D\dT_u,E'_q)$ with $E'\in\Irr_\d\dW$ and summing over all $u\in\dW$ we obtain
$$\align&\sum_i(-1)^i\sum_{w\in\dW}\sum_{\c\in L}\sum_{u\in\dW}q^{l(u)}D(E')(q)\i\\&
\T \tr(\D\dT_u,E'_q)q^{i/2}(\dr'\ot\c:\HH^i(\bX_w))\ti{}\dsg(wu)\dQ_{w,u}(q)\c(\s)\i\\&
=|G^{F'}|\i\sum_{x\in G^{F'}}\tr(\th\i\hx\i,V_{E'})\tr(x\s,\dr')\tag b\endalign$$
for any $\dr'\in\dce_0$, $E'\in\Irr_\d\dW$.

Let $h(r)$ be the common value of the two sides of (b). We now show that (for fixed $\dr'$, $E'$) $h(r)$ has the
following properties:

(i) $|G^{F'}|h(r)$ is a cyclotomic integer;

(ii) when $r$ varies (as a multiple of $dd'$), $h(r)$ is in a fixed (cyclotomic) subfield of $\bbq$ of finite 
degree over $\QQ$;

(iii) all complex conjugates of $h(r)$ have absolute value $\le1$.
\nl
Now (i) is obvious from  the right hand side of (b) since the character values of a finite group (namely 
$\wt{G^F}$ and $\hG^{F'}$) are cyclotomic integers. Also (ii) is obvious from  the left hand side of (b); we use 
that there are only finitely many $\c(\s)$ (roots of $1$) and that $\tr(\D\dT_{u\i},E'_q)\in\QQ(q^{1/2})$. To 
prove (iii) we set:
$$e(r)=|G^{F'}|\i\sum_{x\in G^{F'}}\tr(\th\i\hx\i,V_{E'})\tr(\hx\th,V_{E'})$$
$$e'(r)=|G^{F'}|\i\sum_{x\in G^{F'}}\tr(\s\i x\i,\dr')\tr(x\s,\dr').$$
By the Cauchy-Schwarz inequality, the absolute value square of $h(r)$ (in the form given by the right hand side of
(b)) is $\le$ than $e(r)e'(r)$. Hence it is enough to show that $e(r)=e'(r)=1$. Now $e'(r)=1$ by the 
orthogonality relations for irreducible characters of $\hG^{F'}$ which remain irreducible when restricted to 
$G^{F'}$. It remains to show that $e(r)=1$. Setting $x=aF(a\i),a\in G$ in the definition of $e(r)$ we have 
$$e(r)=|G^{F'}|\i|G^F|\i\sum\Sb a\in G;\\ F'(aF(a\i))=aF(a\i)\endSb
\tr((a\i F'(a)\th)\i,V_{E'})\tr(a\i F'(a)\th,V_{E'}).$$
(We use that, if $a$ is in the sum and $b\in G^F$ then 
$$\tr(a\i F'(a)\th,V_{E'})=\tr(b\i a\i F'(a)F'(b)\th,V_{E'}).$$
 It is enough to show that if $z\in G^F$ then
$\tr(z\th,V_{E'})=\tr(b\i zF'(b)\th,V_{E'})$ which follows from $F'(b)\th=\th b$.) Setting $z=a\i F'(a)\in G^F$ 
we obtain
$$e(r)=|G^F|\i\sum_{z\in G^F}\tr((z\th)\i,V_{E'})\tr(z\th,V_{E'})$$
and this equals $1$ by the orthogonality relations for irreducible characters of $\wt{G^F}$ which remain 
irreducible when restricted to $G^F$. This proves (iii).

If we combine the various imbeddings of the number field in (ii) into $\CC$ we get an imbedding of that field 
into some $\CC^n$. By (i),(iii) the set $\{h(r)\}$ is carried by that imbedding into a discrete and bounded 
subset of $\CC^n$. Hence the set $\{h(r)\}$ is finite. This holds for any $\dr',E'$ which are in finite number. 
We see that there exists an infinite subset $\ci$ of $\{dd',2dd',3dd',\do\}$ such that for any $\dr',E'$, the 
expression (b) is a constant $c_{\dr',E'}$ when $r$ runs through $\ci$. Thus we have
$$\align&\sum_i(-1)^i\sum_{w\in\dW}\sum_{\c\in L}\sum_{u\in\dW}q^{l(u)}D(E')(q)\i\tr(\D\dT_u,E'_q)\\&
\T q^{i/2}(\dr'\ot\c:\HH^i(\bX_w))\ti{}\dsg(wu)\dQ_{w,u}(q)\c(\s)\i=c_{\dr',E'}\tag c\endalign$$
for any $\dr'\in\dce_0$, $E'\in\Irr_\d\dW$, $r\in\ci$.

We multiply the two sides of 2.14(a) by $|G^{F'}|\i\tr(x\s,\dr')$ with $\dr'\in\dce_0$ and sum over all 
$x\in G^{F'}$; we obtain
$$\align&\sum_i(-1)^i\sum_{\c\in L}q^{i/2}(\dr'\ot\c:\HH^i(\bX_w))\ti{}\c(\s)\i\\&
=\sum_{y\in\dW}\dP_{y,w}(q)\sum_{E\in\Irr_\d\dW}\sum_{x\in G^{F'}}
|G^{F'}|\i\tr(x\s,\dr')\tr(\th\i\hx\i,V_E)\tr(\D\dT_y,E_q)\endalign$$
that is 
$$\sum_i(-1)^i\sum_{\c\in L}q^{i/2}(\dr'\ot\c:\HH^i(\bX_w))\ti{}\c(\s)\i
=\sum_{y\in\dW}\dP_{y,w}(q)\sum_{E\in\Irr_\d\dW}c_{\dr',E}\tr(\D\dT_y,E_q)$$
for any $w\in\dW$, $\dr'\in\dce_0$, $r\in\ci$. Since $r$ runs through an infinite set we have an equality of 
polynomials in $v$:
$$\sum_i(-1)^i\sum_{\c\in L}(\dr'\ot\c:\HH^i(\bX_w))\ti{}\c(\s)\i v^i
=\sum_{y\in\dW}\dP_{y,w}(v^2)\sum_{E\in\Irr_\d\dW}c_{\dr',E}\tr(\D\dT_y,E_{v^2})\tag d$$
for any $w\in\dW$, $\dr'\in\dce_0$. Similarly, since (c) holds for $r$ running through an infinite set, we have 
an equality of polynomials in $v$:
$$\align&\sum_i(-1)^i\sum_{w\in\dW}\sum_{\c\in L}\sum_{u\in\dW}v^{2l(u)}D(E')(v^2)\i\\&
\T\tr(\D\dT_u,E'_{v^2})(\dr'\ot\c:\HH^i(\bX_w))\ti{}\dsg(wu)\dQ_{w,u}(v^2)\c(\s)\i v^i=c_{\dr',E'}\endalign$$
for any $\dr'\in\dce_0$, $E'\in\Irr_\d\dW$. Setting $v=1$ in the last equality and using the identity 
$D(E')(1)=|\dW|$ we obtain
$$\align&c_{\dr',E'}=\sum_i(-1)^i\sum_{w\in\dW}\sum_{\c\in L}\sum_{u\in\dW}|\dW|\i\\&
\T\tr(\D u,E')(\dr'\ot\c:\HH^i(\bX_w))\ti{}\dsg(wu)\dQ_{w,u}(1)\c(\s)\i.\tag e\endalign$$
By 2.3(a) we have for fixed $\dr'$:  
$$\sum_i(-1)^i\sum_{\c\in L}(\dr'\ot\c:\HH^i(\bX_w))\ti{}\c(\s)\i=\sum_i(-1)^i(\dr':\HH^i(\bX_w))$$
which by 2.10(a) is equal to
$$\sum_{y\in\dW;y\le w}\dP_{y,w}(1)\sum_i(-1)^i(\dr':H^i_c(X_y)).$$
Substituting this into (e) we obtain
$$\align&c_{\dr',E'}=\sum_{w\in\dW}\sum_{u\in\dW}|\dW|\i\tr(\D u,E')\dsg(wu)\dQ_{w,u}(1)\\&\T|G^{F'}|\i
\sum_{y\in\dW;y\le w}\dP_{y,w}(1)\sum_i(-1)^i(\dr':H^i_c(X_y)).\endalign$$
Using the definition of $\dQ_{w,u}(1)$ we deduce
$$c_{\dr',E'}=\sum_{u\in\dW}|\dW|\i\tr(\D u,E')|G^{F'}|\i(\dr':\dR_u)=(\dr':\dR_{E'})$$
for any $\dr'\in\dce_0$, $E'\in\Irr_\d\dW$. Substituting this into (d) and using the equality
$$\sum_{\c\in L}(\dr'\ot\c:\HH^i(\bX_w))\c(\s)\i=(\dr':\HH^i(\bX_w))$$
(see 2.3(a)), we obtain
$$\align&\sum_i(-1)^i(\dr':\HH^i(\bX_w))v^i
=\sum_{y\in\dW}\dP_{y,w}(v^2)\sum_{E\in\Irr_\d\dW}(\dr':\dR_{E'})\tr(\D\dT_y,E_{v^2})\endalign$$
for any $w\in\dW$, $\dr'\in\dce_0$. The previous equality also holds for $\dr'\in\dce-\dce_0$ (in that case both
sides are zero). This proves 2.4(a) in view of 2.3(b). 

\subhead 2.16\endsubhead
Let $\dot{\le}_{LR}$ be the preorder on $\dW$ (with the weight function $l|_{\dW}$) defined in \cite{\LU, 8.1} 
where it is denoted by $\le_{LR}$ and let $\dot{\si}_{LR}$ be the corresponding equivalence relation on $\dW$ 
(the equivalence classes are the two-sided cells of $\dW$ with the weight function above; they form a set 
$\dcc$). For $w,w'\in\dW$ we write $w'\dot{<}_{LR}w$ when $w'\dot{\le}_{LR}w$ and $w'\not\dot{\si}_{LR}w$. The 
relations $\dot{\le}_{LR},\dot{<}_{LR}$ on $\dW$ induce relations on $\dcc$ denoted again by $\dot{\le}_{LR}$, 
$\dot{<}_{LR}$. It is known \cite{\LU, \S10} that if $\fc\in\dcc$ then $\fc^*:=\fc s_I=s_I\fc\in\dcc$ and that 
$\fc\m\fc^*$ reverses the preorder $\dot{\le}_{LR}$. If $E\in\Irr\dW$ then there is a unique $\fc\in\dcc$ such 
that $\dt_w:E_\sp@>>>E_\sp$ is nonzero for some $w\in\fc$ and is zero for all $w\in\dW-\fc$; we then write 
$\fc=\fc_E$. 

Let $a:\dW@>>>\NN$ is the function defined as in \cite{\LU, 13.6} in terms of the weight function $l|_{\dW}$; 
this function is in fact the restriction to $\dW$ of the function $a:W@>>>\NN$ in 1.2, see \cite{\LU, 16.5}. Let 
$E\m a_E$ be the function $\Irr\dW@>>>\NN$  defined in \cite{\LU, 20.6} (in terms of the weight function of
$\dW$). For $E\in\Irr\dW$ we set $a'_E=a_{E^\da}$.

Using \cite{\LU, \S20} we see that for $w\in\dW, E\in\Irr_\d\dW$ we have
$$v^{-l(w)}\tr(\D\dT_w,E_{v^2})=\dc'_{w,\D,E}v^{a'_E}+\text{lower powers of $v$}$$
where $\dc'_{w,\D,E}\in\ZZ$; by an argument similar to that in \cite{\LU, 20.10}, we have 
$\dc'_{w,\D,E}=\tr(\D',E^\da_\sp)$ where $\D':E^\da@>>>E^\da$ is $\D\ot1:E\ot\dsg@>>>E\ot\dsg$.

Let $\dcar_{\dW}$ be the vector space of formal linear combinations $\sum_{E\in\Irr_\d\dW}r_EE$, $r_E\in\QQ$. For 
$w\in\dW$ we set
$$\dfA_w=\sum_{E\in\Irr_\d\dW}\dc'_{w,\D,E}E\in\dcar_{\dW}.$$
From the definitions we see that

(a) {\it If $\fc\in\dcc,w\in\fc$ then $\dfA_w$ is a linear combination of $E\in\Irr_\d\dW$ such that 
$\fc_E=\fc^*$.}
\nl
Conversely, if $E\in\Irr_\d\dW$ then

(b) {\it $E$ is a $\QQ$-linear combination of elements $\dfA_x$ ($x\in\fc_E^*$).}
\nl
The proof is along the lines of that of \cite{\ORA, 5.13(ii)} (which is the analogous result for $W$ instead of 
$\dW$).

For any $\x=\sum_{E\in\Irr_\d\dW}r_EE\in\dcar_{\dW}$ we set $\dR_\x=\sum_{E\in\Irr_\d\dW}r_E\dR_E\in[\dce]$. In 
particular, for $w\in\dW$, $\dR_{\dfA_w}\in[\dce]$ is defined.

\proclaim{Proposition 2.17} Let $w\in\dW$. We have
$$\align&(-1)^{l(w)-a(w)}\HH^{l(w)-a(w)}(\bX_w)=(-1)^{l(w)+a(w)}\HH^{l(w)+a(w)}(\bX_w)\\&=
\dR_{\dfA_w}+\sum_{w'\in\dW,w'\dot{<}_{LR}w}n_{w',w}\dR_{\dfA_{w'}}\in[\dce]\tag a\endalign$$
for certain rational numbers $n_{w',w}$.
\endproclaim
The first equality follows from the Lefschetz hard theorem in intersection cohomology applied to $\bX_w$ which
has pure dimension $l(w)$. In the rest of the proof we concentrate on the second equality. Using the second part 
of 2.4 with $i=l(w)+a(w)$ we see that it is enough to show that
$$\align&\sum_{E\in\Irr_\d\dW}\tr(\D\sum_{y\in\dW}\dP_{y,w}(v^2)\dT_y,E_{v^2};l(w)+a(w))E\\&=
\dfA_w+\sum_{w'\in\dW,w'\dot{<}_{LR}w}\dfA_{w'}\tag b\endalign$$
(equality in $\dcar_{\dW}$). 
Let $\fc\in\dcc$ be such that $w\in\fc$. From the definitions we see that, for $E\in\Irr_\d\dW$, the operator 
$\sum_{y\in\dW}\dP_{y,w}(v^2)\dT_y:E_{v^2}@>>>E_{v^2}$ is zero unless $\fc_E^*\dot{\le}_{LR}\fc$. Thus the sum in
the left hand side of (b) can be restricted to the $E$ such that $\fc_E^*\dot{\le}_{LR}\fc$. Next we show that 
for any $E\in\Irr_\d\dW$ we have
$$\tr(\D\sum_{y\in\dW}\dP_{y,w}(v^2)\dT_y,E_{v^2})=\dc'_{w,\D,E}v^{l(w)+a'_E}+\text{lower powers of $v$}.\tag b$$
(Using the definition of $\dc'_{w,\D,E}$, it is enough to show that for any $y\in\dW$, $y<w$ we have
$$\dP_{y,w}(v^2)\tr(\D\dT_y,E_{v^2})\in v^{l(w)+a'_E-1}\ZZ[v\i];$$
since 
$$\tr(\D\dT_y,E_{v^2})\in v^{l(y)+a'_E}\ZZ[v\i],$$
it is enough to note that $\dP_{y,w}(v^2)\in v^{l(w)-l(y)-1}\ZZ[v\i]$.) It follows that any $E$ which appears 
with nonzero coefficient in the first sum in (b) satisfies $\fc_E^*\dot{\le}_{LR}\fc$ and that the contribution 
to that sum of the $E$ such $\fc_E^*=\fc$ (hence $a'_E=a(w)$, see \cite{\LU, 20.6(c)}) is 
$$\sum_{E\in\Irr_\d\dW;\fc_E^*=\fc}\dc'_{w,\D,E}E=\dfA_w$$
(see 2.16(a)). It remains to show that if $E$ satisfies $\fc_E^*\dot{<}_{LR}\fc$ then $E$ is a linear combination
of elements $\dfA_{w'}$ with $w'\in\dW,w'\dot{<}_{LR}w$. This follows from 2.16(b).

\subhead 2.18\endsubhead
For any $c\in\cc$ we set 
$$[\dce]^{\le c}=\op_{c'\in\cc;c'\le_{LR}c}[\dce]^{c'},\qua[\dce]^{\ge c}=\op_{c'\in\cc;c\le_{LR}c'}[\dce]^{c'}.$$
By \cite{\LU, 23.5} any $\fc\in\dcc$ is contained in a unique two-sided cell $\fc^!\in\cc$; moreover, the map
$\dcc@>>>\cc$, $\fc\m\fc^!$ is injective (we have $\fc^!\cap\dW=c$). Also $(\fc^*)^!=(\fc^!)^*$. (Indeed, if 
$w\in\fc$ so that $w\in\fc^!$, we have $ws_I\in\fc^*$ so that $ws_I\in\fc^*$ and $ws_I\in(\fc^*)^!$; but we have 
also $ws_I\in\fc^!s_I=(\fc^!)^*$. Hence the desired equality.)

We show:

(a) If $\fc\in\dcc$ and $w\in\fc$ then $\dR_{\dfA_w}\in[\dce]^{\le\fc^!}$.
\nl
We can assume that this is true if $w$ is replaced by $w'$ with $w'\in\dW$, $w'\dot{<}_{LR}w$. By 1.7, the 
$G^{F'}$-module $\HH^{l(w)-a(w)}(\bX_{w})$ is a sum of representations isomorphic to $\r$ such that 
$\un{\r}^*\le_{LR}\fc^!$. Hence $\HH^{l(w)-a(w)}(\bX_{w})\in[\dce]^{\le\fc^!}$. Using now 2.17 we see that 

(b) $\dR_{\dfA_w}+\sum_{w'\in\dW,w'\dot{<}_{LR}w}n_{w',w}\dR_{\dfA_{w'}}\in[\dce]^{\le\fc^!}$
\nl
with $n_{w',w}\in\QQ$. By the induction hypothesis for any $w'$ in the last sum (with $w'\in\fc',\fc'\in\dcc$) we 
have $\dR_{\dfA_{w'}}\in[\dce]^{\le\fc'{}^!}$. By arguments in \cite{\LU, \S16} (see especially 16.6, 16.13(a)), 
from $w'\dot{\le}_{LR}w$ we deduce that $w'\le_{LR}w$ that is $\fc'{}^!\le\fc^!$. Hence 
$[\dce]^{\le\fc'{}^!}\sub[\dce]^{\le\fc^!}$. Thus $\dR_{\dfA_{w'}}\in[\dce]^{\le\fc^!}$. Introducing this in (b) 
we see that  $\dR_{\dfA_w}\in[\dce]^{\le\fc^!}$. This proves (a).

We show:

(c) {\it If $E\in\Irr_\d\dW$ then $\dR_E\in[\dce]^{\le(\fc_E^*)^!}$.}
\nl
This follows from (a) using 2.16(b).

\subhead 2.19\endsubhead
Recall the operation of "duality" \cite{\DLL,\DLLL} $D:[\ce]@>>>[\ce]$ which is a linear map such that
$D(\r)=\ti{\r}$ where $\r\m\ti{\r}$ is an involution of $\ce$ and such that $D(R_w)=\sg(w)R_w$ for any $w\in W$.
It follows that if $E\in\Irr_\d W$ then $D(R_E)=\pm R_{E^\da}$. Hence if $\r\in\ce$, then 
$\ti{\un{\r}}=\un{\r}^*$.

The operator $D$ is generalized in \cite{\DM,\S3} to a linear map $\dD:[\dce]@>>>[\dce]$ with the following
properties: (i) if $\dr\in\dce_0$ and $\dr|_{G^{F'}}=\r\in\ce_0$ then $\dD(\dr)$ is up to scalar of the form
$\dr'\in\dce_0$ where $\dr'|_{G^{F'}}=\ti{\r}$ (as above); (ii) $\dD(\dR_w)=\dsg(w)\dR_w$ for any $w\in\dW$. It 
follows that if $E\in\Irr_\d\dW$ then $\dD(\dR_E)=\pm\dR_{E^\da}$.
 
We show:

(a) {\it if $c\in\cc$ then $\dD([\dce]^{\le c})=[\dce]^{\ge c^*}$.}
\nl
It is enough to show that, if $\r\in\ce_0,\ti{\r}\in\ce_0$ are as above, then $\un{\r}^*\le_{LR}c$ implies
$\ti{\un{\r}}^*\ge_{LR}c^*$, that is $\un{\r}\ge_{LR}c^*$. But this follows from the fact ${}^*$ reverses the 
preorder 
$\le_{LR}$.

Applying $\dD$ to 2.18(c) with $E$ replaced by $E'$ (as above) and using (a) we obtain
$\dR_E\in[\dce]^{\ge((\fc_{E'}^*)^!)^*}$. We have $((\fc_{E'}^*)^!)^*=(\fc_E^!)^*=(\fc_E^*)^!$. Hence

(b) $\dR_E\in[\dce]^{\ge(\fc_E^*)^!}$. 
\nl
From (b) and 2.18(b) we deduce that $\dR_E\in[\dce]^{\le(\fc_E^*)^!}\cap[\dce]^{\ge(\fc_E^*)^!}$. Hence

(c) $\dR_E\in[\dce]^{(\fc_E^*)^!}$.
\nl
This proves 2.4(ii).

From (c) we deduce, using 2.16(a), that for any $\fc\in\dcc$ and any $w\in\fc$ we have

(d) $\dR_{\dfA_w}\in[\dce]^{\fc^!}$.

\subhead 2.20\endsubhead
Let $\fc\in\dcc,w\in\fc$. Let $\p_{\fc^!}:[\dce]@>>>[\dce]^{\fc^!}$ be the canonical projection. We apply 
$\p_{\fc^!}$ to the  equality between the first and third member of 2.17(a) and we use 2.19(d). We obtain
$$(-1)^{l(w)-a(w)}\p_{\fc^!}(\HH^{l(w)-a(w)}(\bX_w))=\dR_{\dfA_w}.$$
We have used that if $w'\in\dW,w'\dot{<}_{LR}w$ then $\p_{\fc^!}(\dR_{\dfA_{w'}})=0$. (Indeed, let $\fc'\in\dcc$
be such that $w'\in\fc'$. We have $\fc'\ne\fc$. Using 2.19(d) for $w'$ instead of $w$ it is enough to show that 
$\fc'{}^!\ne\fc^!$; this follows from the injectivity of the map $\fc\m\fc^!$.) Note that 
$\HH^{l(w)-a(w)}(\bX_w)_\di$ (see 1.6) is a $\hG^{F'}$-submodule of $\HH^{l(w)-a(w)}(\bX_w)$; moreover from the 
definitions we have $\p_{\fc^!}(\HH^{l(w)-a(w)}(\bX_w))=\HH^{l(w)-a(w)}(\bX_w)_\di$ (equality in $[\dce]$). Thus 
we have the following generalization of 1.6(c):
$$(-1)^{l(w)-a(w)}\HH^{l(w)-a(w)}(\bX_w)_\di=\dR_{\dfA_w}.\tag a$$
In particular, the following generalization of 1.6(d) holds:
$$(-1)^{l(w)-a(w)}\dR_{\dfA_w}\text{ can be represented by an actual $\hG^{F'}$-module}.\tag b$$

\subhead 2.21\endsubhead
Assume now that $G$ is an even special orthogonal group over $F_{q'}$ and that $\s$ is conjugation by a reflection
defined over $F_{q'}$ so that $\hG$ is an even full orthogonal group over $F_{q'}$. Using the results of this 
paper (especially 2.20(b)) one can show that if $E\in\Irr_\d\dW$ then there exists $c\in\cc$ such that 
$$2^n\dR_E=\sum_{\r\in\ce;\un{\r}=c}\ti{\r}$$
where $n$ is defined by $|\{\r\in\ce;\un{\r}=c\}|=2^{2n}$ and for each $\r$ in the sum, $\ti{\r}$ denotes a 
certain extension of $\r$ to a $G^{F'}$-module. The proof will be given elsewhere.

\enddocument